\documentclass[12pt]{article}
\usepackage[centertags]{amsmath}
\usepackage{amsfonts}
\usepackage{amssymb}
\usepackage{amsthm}
\input{amssym.def}

\newtheorem{Definition}{Definition}[section]
\newtheorem{Theorem}[Definition]{Theorem}
\newtheorem{Lemma}[Definition]{Lemma}

\newtheorem{Corollary}[Definition]{Corollary}

\newtheorem{Example}[Definition]{Example}

\newcommand{\lc}{\mathcal{L}}
\newcommand{\rc}{\mathcal{R}}
\newcommand{\hc}{\mathcal{H}}
\newcommand{\jc}{\mathcal{J}}

\title{\Large \bf On inverse and right inverse ordered semigroups}
\author{A. Jamadar and K. Hansda}
\date{}

\begin{document}

\maketitle

\begin{abstract}{\footnotesize}
A regular ordered semigroup $S$ is called right inverse if every
principal left ideal of $S$ is  generated by an $\rc$-unique ordered
idempotent. Here we explore the theory of right inverse ordered
semigroups. We show that a regular ordered semigroup is right
inverse if and only if any two right inverses of an element $a\in S$
are $\rc$-related. Furthermore, different characterizations of right
Clifford, right group-like, group like ordered semigroups are done
by right inverse ordered semigroups. Thus a foundation of right
inverse semigroups has been developed.

\end{abstract}

\section{Introduction}

Right inverse semigroups are those, every element of which has
unique right inverse. Thus naturally it becomes generalization of
inverse semigroups. Many extensive studies have been done on right
inverse semigroups by P.S. Venkatesan \cite{VEN}, G.L. Bailes
\cite{ba} and some others. P.S. Venkatesan \cite{VEN} studied these
semigroups under the name of right unipotent semigroups. He showed
that a semigroup is right inverse if and only if every right ideal
of it generated by an idempotent.

T. Saito \cite{Saito 1971}, studied inverse semigroup by introducing
simple ordered on it. Bhuniya and Hansda  \cite{bh1} have deal with
ordered semigroups in which any two inverses of an element are
$\hc$-related. These ordered semigroups are the analouge of inverse
semigroups. Hansda and Jamadar \cite{HJ} named these ordered
semigroups inverse ordered semigroups. They gave a detailed
exposition on the characterization of these ordered semigroups. Here
we generalize such ordered semigroups into right inverse ordered
semigroups. This paper is inspired by the works done by
P.S.Venkatesan \cite{VEN}, G.L.Bailes \cite{ba}.

The presentation of this article is as follows: This section is
followed by preliminaries. Section 3 is devoted to the right inverse
ordered semigroups. Here Clifford ordered semigroups have been
characterized by right inverse semigroups.

\section{Preliminaries}
 An ordered semigroup  is a partiality ordered set
$(S,\cdot,\leq)$, and at the same time a semigroup $(S,\cdot)$ such
that for all $a, \;b , \;x \in S  \;a \leq b$ implies
 $xa\leq xb \;\textrm{and} \;a x \leq b x $.
It is denoted by $(S,\cdot, \leq)$. For every subset $H\subseteq S$,
denote $(H]=\{t\in S: t\leq h, \;\textrm{for some} \;h\in H\}$.
Throughout this article, unless stated otherwise, $S$ stands for an
ordered semigroup and We assume that $S$ does not contain the zero
element.

An equivalence relation $\rho$ is called left (right) congruence if
for $a, b, c \in S \;a\rho b \;\textrm{ implies} \;ca \rho cb \;(ac
\rho bc)$.  By a congruence we mean both left and right congruence.
A congruence $\rho$ is called semilattice congruence on $S$ if  for
all $a, b \in S, \;a \rho a^{2} \;\textrm {and} \;ab \rho \;ba$. By
a complete semilattice congruence we mean a semilattice congruence
$\sigma$ on $S$ such that for $a, b \in S, \;a \leq b$ implies that
$a \sigma ab$. The ordered semigroup $S$ is called complete
semilattice of subsemigroups of type $\tau$ if there exists a
complete semilattice congruence $\rho $ such that $(x)_{\rho}$ is a
type $\tau$ subsemigroup of $S$. Let $I$ be a nonempty subset of
$S$. Then $I$ is called a left(right) ideal of $S$, if $SI\subseteq
I (IS\subseteq I)$ and $(I]\subseteq I$. If $I$ is both  left and
right ideal, then it is called an ideal of $S$. We call $S$ a (left,
right) simple ordered semigroup if it does not contain any proper
(left,right) ideal.  For $a\in S$, the smallest (left, right) ideal
of $S$ that contains $a$ is denoted by $(L(a),  \;R(a)) I(a)$.

$S$ is said to be regular (resp. Completely regular, right regular)
ordered semigroup if for every $a \in S, \;a\in (aSa](a\in
(a^2Sa^2], \;a\in (a^2S])$. Due to Kehayopulu \cite{Ke2006} Green's
relations on a regular ordered semigroup given as follows:

$a\lc b\; \textrm{if} \;L(a)= L(b)$, $a\rc b \;\textrm{if} \;R(a)=
R(b)$, $a\jc b \;\textrm{if} \;I(a)= I(b)$, $\hc=\lc \cap \rc$.

This four relation $\lc, \rc, \jc,\textrm{and}  \;\hc$ are
equivalence relation.

A regular ordered semigroup $S$ is said to be group-like (resp. left
group-like) \cite{bh1} ordered semigroup if for every $a,b\in S,
\;a\in (Sb] \;\textrm{and} \;b \in (aS](\textrm{resp}. \;a \in
(Sb])$. Right group like ordered semigroup can be defined dually. A
regular ordered semigroup $S$ is called a right (left) Clifford
 \cite{bh1} ordered semigroup if for all $a\in S, (Sa]\subseteq
(aS],( (aS]\subseteq (Sa])$. Every right (left) group like ordered
semigroup is a right (left) Clifford ordered semigroup. An element
$b\in S$ is said to be an inverse of $a\in S$ if $a\leq aba$ and
$b\leq bab$. The set of all inverses of an element $a$ is denoted by
$V_\leq (a)$.
\begin{Theorem} \cite{bh1}
Let S be a regular ordered semigroup. Then the following statements
are equivalent.
\begin{enumerate}
\item \vspace{-.4cm}
$S$ is right Clifford ordered semigroup;
\item \vspace {-.4cm}
for all $e\in E_\leq(S), (Se]\subseteq (eS]$;
\item \vspace{-.4cm}
for all $a\in S$, and $e\in E_\leq(S)$, there is $x\in S$ such that
$ea\leq ax$;
\item \vspace{-.4cm}
for all $a,b \in S$, there is $x\in S$ such that $ba\leq ax$;
\item \vspace{-.4cm}
$\lc\subseteq \rc$ on $S$.
\end{enumerate}
\end{Theorem}
\begin{Lemma}\label{6}\cite{bh1}
Let $S$ be a right Clifford ordered semigroup. Then the following
conditions hold in $S$.
\begin{enumerate}
\item \vspace{-.4cm}
$a\in (a^2Sa]$, for every $a\in S$;
\item \vspace{-.4cm}
$ef\in (feSef]$, for every $e,f\in E_\leq (S)$.
\end{enumerate}
\end{Lemma}
\begin{Theorem} \cite{bh1}\label{1}
Let $S$ be an ordered ordered semigroup. Then $S$ is right (left)
Clifford ordered semigroup if and only if $\rc (\lc)$ is the least
complete semilattice congruence on $S$.
\end{Theorem}
\begin{Theorem} \cite{bh1}\label{2}
Let $S$ be a regular ordered semigroup. Then $S$ is right (left)
Clifford ordered semigroup if and only if it is a complete
semilattice of right (left) group like ordered semigroups.
\end{Theorem}

\section{Right Inverse ordered semigroup}
\subsection{Right inverse ordered semigroups}
Let $S$ be an ordered semigroup and $\rho$ be an equivalence
relation on $S$. In broad sense $\rho$-unique we shall mean the
uniqueness in respect of the relation $\rho$. For example consider a
subset $T$ of $S$ such that $a,b$ are  generators of $T$. Now if
$\;a\rho b$ we say that $T$ is generated by $\rho$-unique element
$a$.
\begin{Definition}
A regular ordered semigroup $S$ is called right inverse if every
principal left ideal is generated by an $\rc-$unique ordered
idempotent of $S$.
\end{Definition}

We now present results on the  role of ordered idempotents to
characterize right inverse ordered semigroups.
\begin{Theorem}
A regular ordered semigroup $S$ is a right inverse if and only if
for any two idempotents $e,f\in E_\leq(S)$, $e\lc f$ implies $e\hc
f$.
\end{Theorem}
\begin{Theorem}
Let $S$ be a regular ordered semigroup. Then $S$ is left (right)
group like ordered semigroup if and only if any two ordered
idempotents are $\lc \;(\rc)-$related.
\end{Theorem}

\begin{Corollary}
Every right inverse left group like ordered semigroup is a group
like ordered semigroup.
\end{Corollary}
\begin{Theorem}
Let $S$ be a regular ordered semigroup. Then any two inverses of an
element are $\lc-$ related if and only if $ef\in(eSfSe]$; for some
$e,f\in E_{\leq}(S)$.
\end{Theorem}

In the following theorem, we have shown that any two inverses of an
element are $\rc$-related in a right inverse ordered semigroup. So
in the broad sense they are $\rc$-unique.

\begin{Theorem}\label{5}
The following conditions are equivalent on a regular ordered
semigroup $S$.
\begin{enumerate}
\item\vspace{-.4cm}
$S$ is right inverse;
\item\vspace{-.4cm}
for $a\in S$ and $a', a'' \in V_{\leq}(a)$, $a'\rc a''$;
\item\vspace{-.4cm}
for $e,f\in E_{\leq}(S)$, $ef\in (fSeSf]$;
\item\vspace{-.4cm}
$(eS] \cap (fS]=(efS]$;
\item\vspace{-.4cm}
for $e\in E_{\leq}(S)$ and $x\in (Se]$\textrm{ implies} $x'\in
(eS]$, where $x\in S$ and $x'\in V_{\leq}(x)$.
\end{enumerate}
\end{Theorem}

\begin{Corollary}
Let $S$ be a right inverse ordered semigroup. Then any two ordered
idempotents $e,f\in E_\leq(S)$ are $\hc$- commutative if and only if
$(Se]\cap (Sf]= (Sef]$.
\end{Corollary}

\begin{Example}
The ordered semigroup $S= \{a, e, f  \}$ defined by multiplication
and order below.
\begin{center}
\begin{tabular}{|l|l|l|l|}
  \hline
  $\cdot$ & $a $&$e$ & $f$ \\
  \hline
  $a$ & $a$ & $e$ & $f$ \\
  \hline
  $e$ & $a$ & $e$ & $f$ \\
  \hline
  $f$ & $a$ & $e$ & $f$ \\
  \hline
\end{tabular}
\end{center}
$$'\leq ' := \{(a,a), (a,e),  (a,f),  (e,e),  (f,f)\}.$$

From above table it is clear that $a, e, f\in E_\leq (S)$. Here
$ae\leq aee= eaee= eaaee$. So $ae\in (eSaSe]$. Also $a\leq aa$
implies that $ea\leq eaa= aeaa= aeeaa$. So $ea\in (aSeSa]$.
Similarly $af\in (fSaSf]$ and $fa\in (aSfSa]$. Also $ef= fef= feeff$
that is $ef\in (fSeSf]$. Similarly it can be shown that $fe\in
(eSfSe]$. Thus $(S, \cdot, \leq)$ is a right inverse ordered
semigroup.
\end{Example}
\begin{Theorem}
Let $F$ be a semigroup. Then the ordered semigroup $P_f(F)$ of all
subsets of $F$ is a right inverse ordered semigroup if and only if
$F$ is a right inverse semigroup.
\end{Theorem}

\begin{Theorem}
Let $S$ be a regular ordered semigroup. Then $S$ is a right inverse
ordered semigroup if and only if $L_e \subseteq (R_e)'$ for any
idempotent $e$ in $S$.
\end{Theorem}

\begin{Theorem}
An ordered semigroup $S$ is right Clifford if and only if $S$ is
right inverse and for every $a\in S$, $a\in (a^2Sa]$.
\end{Theorem}
\begin{Theorem}\label{3}
Let $S$ be a right inverse ordered semigroup. If $S$ is left
Clifford then $S$ is union of group like ordered semigroups.
\end{Theorem}

In the following we show that in a right inverse ordered semigroup
$\rc$ is a congruence if and only if $\lc= \hc$.
\begin{Theorem}\label{4}
Let $S$ be a right inverse ordered semigroup. The following are
equivalent:
\begin{enumerate}
\item\vspace{-.4cm}
$\rc$ is a congruence on $S$;
\item\vspace{-.4cm}
$\lc= \hc$;
\item\vspace{-.4cm}
$S$ is a complete semilattice of right group like ordered
semigroups.
\end{enumerate}
\end{Theorem}

Our paper ends up with the corollary that follows from  Theorem
$\ref{4}$ and Theorem $\ref{3}$ and which gives a characterization
on right inverse semigroup to become a completely regular ordered
semigroup.
\begin{Corollary}
Let $S$ be a right inverse and left regular ordered semigroup. Then
following conditions are equivalent.
\begin{enumerate}
\item\vspace{-.4cm}
$\rc$ is a congruence on $S$;
\item\vspace{-.4cm}
$\lc= \hc$;
\item\vspace{-.4cm}
$S$ is a complete semilattice of right group like ordered
semigroups;
\item\vspace{-.4cm}
$S$ is completely regular.
\end{enumerate}
\end{Corollary}

\bibliographystyle{plain}

\begin{thebibliography}{10}
\baselineskip 5mm



\bibitem{bh1}
A. K. Bhuniya and K. Hansda, {Complete semilattice of ordered
semigroups}, Communicated.

\bibitem{ba}
G. L. Bailes, {Right inverse Semigroups}, \emph{Journal of Algebra},
\textbf{26}(1973), 492-507.


\bibitem{ Howie 1995}
J.M. Howie, Fundamentals of Semigroup Theory, \emph{Clarendon Press,
Oxford},  \textbf 1995.

\bibitem{Ke2 90}
N.Kehayopulu, {Remark in ordered semigroups}, \emph{Math. Japonica}
 \textbf{35}(1990), 1061-1063.


 \bibitem{Ke2002}
N.Kehayopulu and M.Tsingelis, {On Left Regular Ordered Semigroups},
\emph{Southeast Asian Bulletin of Mathematics}
 \textbf{25}(2002),609-615.


\bibitem{Ke2006}
N.Kehayopulu, {Ideals and Green's relations in ordered semigroups},
\emph{International Journal of Mathematics and Mathematical Sciences
}, \textbf{}(2006), 1-8,Article ID 61286.



\bibitem{ke1}
N. Kehayopulu and   Tsingelis,  {Semilattices of Archimedean ordered
Semigroups}, \emph{Algera Colloquium} \textbf{15:3}(2008), 527-540


\bibitem{Ke2009}
 N.Kehayopulu, {Archimeadean ordered semigroups as ideal extensions},\emph{Semigroup Forum
}, \textbf{78}(2009), 343-348.

\bibitem{HJ}
K. Hansda and A. Jamadar, {On inverse semigroups}, Communicated.


\bibitem{Petrich 1977}
M.Petrich, {Lectures in Semigroups}, \emph{John Wiley and Sons },
 1977 .


\bibitem{Saito 1962}
T.Saito, {Ordered idempotent semigroups}, \emph{J. Math. Soc.
Japan}, \textbf{14(2)}(1962), 150-169.

\bibitem{Saito 1971}
T.Saito, {Ordered inverse semigroups}, \emph{Trans. Amer. Math.
Soc}, \textbf{153}(1971),99-138.

\bibitem{VEN}
P. S. Venkatesan, {Right(Left)Inverse Semigroups}, \emph{Journal of
Algebra}, \textbf{31}(1974), 209-217.

\bibitem{VEN1}
P. S. Venkatesan, {On Right Unipotent Semigroups}, \emph{Pacific
Journal of Mathematics}, \textbf{63}(1976), 555-561.

\bibitem{VEN2}
P. S. Venkatesan, {On Right Unipotent Semigroups II}, \emph{Glasgow
Math.J.}, \textbf{19}(1978), 63-68.

\bibitem{VEN3}
P. S. Venkatesan, {Bisimple Left Inverse Semigroups},
\emph{Semigroup Forum}, \textbf{4}(1972),34-45.


\end{thebibliography}

\end{document}